\documentclass[12pt, twoside, leqno]{article}
\usepackage{amsmath,amsthm}
\usepackage{amssymb,latexsym}
\usepackage{enumerate}
\usepackage{float}
\usepackage{epsfig}
\usepackage{color}
\usepackage{graphics,graphicx}
\usepackage{boxedminipage}

\newtheorem{thm}{Theorem}[section]

\newtheorem{lem}[thm]{Lemma}

\newtheorem{prop}[thm]{Proposition}

\newtheorem{mainthm}[thm]{Main Theorem}

\theoremstyle{definition}
\newtheorem{defi}[thm]{Definition}

\numberwithin{equation}{section}

\newcommand{\R}{{\bf R}}

\newcommand{\La}{\Lambda}
\newcommand{\N}{{\rm I\kern-4pt N}}
\newcommand{\RR}{{\rm I\kern-4pt R}}
\newcommand{\NN}{{\rm I\kern-4pt N}}

\newcommand{\al}{\alpha}

\newcommand{\Ome}{\Omega}

\newcommand{\ds}{\displaystyle}

\def \ind {\hbox{ 1\hskip -3pt I}}
\newcommand{\BE}{\begin{equation}}
\newcommand{\EE}{\end{equation}}

\begin{document}


\baselineskip=17pt


\title{Analysis of a fractal boundary: the graph of the Knopp function}

\author{Mourad Ben Slimane\\
Department of Mathematics, College of Science,
King Saud University,\\ P. O. Box 2455, Riyadh 11451, Saudi Arabia,\\ Email:  mbenslimane@ksu.edu.sa\\
\and
Clothilde M\'elot\\
Aix-Marseille Univ,  CNRS, LATP, UMR 6632,\\
  F-13453 Marseille, France,\\
   Email: melot@cmi.univ-mrs.fr}

\date{}

\maketitle


\renewcommand{\thefootnote}{}



\renewcommand{\thefootnote}{\arabic{footnote}}
\setcounter{footnote}{0}


\begin{abstract}
A usual classification tool to study a fractal interface is the computation of its fractal dimension.
 But a recent method developed by Y. Heurteaux and S. Jaffard proposes to compute either weak and strong accessibility exponents
 or local $L^p$ regularity exponents (the so-called $p$-exponent). These exponents describe locally the behavior of the interface.
We apply this method to the graph of the Knopp function which is defined for $x\in[0,1]$  as
$F(x) = \sum_{j=0}^\infty  2^{-\alpha j}  \phi(2^{j}x)$
where $0<\alpha <1$  and $\phi(x) = dist(x,\mathbb{Z})$.  The Knopp function itself has everywhere the
same $p$-exponent $\alpha$. Nevertheless,
using the characterization of the maxima and minima done by B. Dubuc and S. Dubuc,
we will compute the $p$-exponent  of the characteristic function  of domain
under the graph of $F$ at each point $(x,F(x))$ and show that $p$-exponents, weak and strong accessibility exponents change from point to point. Furthermore we will derive a characterization of the local extrema of the function according to the values of these exponents.
\end{abstract}

%
%
%
%
%
%
%
%
%

 \section{Introduction}
At the beginning of the century several examples of non differentiable functions were studied, such as the Weiertrass function or the example we will focus on in the following, i.e the Takagi-Knopp or so called Knopp function (see \cite{All-Kaw} and references therein for a review). The issue was the study of the regularity.\\


Indeed in 1918, Knopp \cite{Knopp} introduced  a new family of non differentiable functions defined on the interval $[0,1]$.
Going beyond the construction of Weierstrass of a continuous non differentiable function, his goal was to build examples of
continuous functions  for which one sided limits of the difference quotient at all points don't exist. He considered  the function
  $F_{a,b}$ given by the series, for $x\in [0,1]$
\BE \label{expansion0}
F_{a,b}(x) = \sum_{j=0}^\infty  a^{ j}  \phi(b^{j}x)
\EE
where $\phi(x) = dist(x,\mathbb{Z})$, $0<a<1$, $b$ is an integer such that $ab>4 $.\\

 For $b=2$ and $a=2^{-\alpha}$, this function can be seen as a series
  expanded in the Faber-Schauder basis
   $\La_{j,k}:x\mapsto 2^{\frac{j}{2}}\La(2^{j}x-k), j \in \NN, k=0,\cdots,2^{j}-1$,
   where $\La$ is the Schauder function defined by  $\La(x)=\inf(x,1-x)$ if
$x \in [0,1]$ and  0 elsewhere. In fact
\BE \label{expansion0}
F_{2^{-\alpha},2}(x) =
\sum_{j=0}^\infty \sum_{k=0}^{2^{j}-1} 2^{ -\alpha j}  \La(2^{j}x-k)\;.
\EE
  {\it We will write $F$ for $ F_{2^{-\alpha},2}$ in the following. }\\
  
   Thus, for example, using the characterization of Lipschitz
  spaces with the help of coefficients in the Schauder-basis  \cite{Cies}, one gets immediately the fact that $F$ belongs to $C^\alpha([0,1])$.\\
  
A further step to study the regularity of this function can be to follow the ideas developped in multifractal analysis. The goal in multifractal analysis is to study the sets of points where the function has a given pointwise regularity, and doing so checking if the regularity changes from point to point and quantify these changes. Recall the definition of H\"older pointwise regularity and local $L^p$ regularity.

\begin{defi}
Let $x_0 \in \R^d$ and $\alpha\geq 0$.
A locally bounded function $f: \R^d \rightarrow \R$
 belongs to $C^{\alpha}(x_0)$ if there exists $C >0$ and a polynomial $P=P_{x_0}$ with $deg(P)\leq [{\alpha}]$,
  such that on a neighborhood of $x_0$,
\begin{equation}
\label{Hol}  |f(x)-P_{x_0}(x)| \leq C |x-x_0|^{\alpha}.
\end{equation}
 The pointwise H\"older exponent of $f$ at $x_0$ is $ h_f(x_0)=\sup\{\alpha:f\in C^{\alpha}(x_0)\}$.
\end{defi}

\begin{defi}\label{defi:pexponent}\cite{CZ61}
Let $x_0\in\R^d$.
 Let $p\in[1,\infty]$ and $u$ such that $u\geq -\frac{d}{p}$.
Let $f$ be a function in $L^p_{loc}$.
The function $f$ belongs to $T^p_u(x_0)$ if there exists $R>0$,  a polynomial $P$ with $deg(P)\leq u$, and $C>0$ such that
\begin{equation}\label{tp}
\forall \rho\leq R:\mbox{ }\left(\frac{1}{\rho^d}\int_{|x-x_0|\leq
\rho}|f(x)-P(x)|^p
dx\right)^{\frac{1}{p}}\leq C\rho^{u}.
\end{equation}

The $p$-exponent of $f$ at $x_0$ is $u^p_f(x_0)=\sup\{u:f\in T^p_u(x_0)\}$.

\end{defi}

Then again with the help of the Faber-Schauder basis one can prove that for all $x_0\in[0,1]$, $F$ is in $C^\alpha(x_0)$ ( details for this technique can be found in \cite{JaffMand}). It is then easy to check that actually $u_p^f(x_0)=h_f(x_0)=\alpha $ at all $x_0\in[0,1]$. 
Thus from the point of view of various notions of regularity, even if it is not differentiable,
the function $F$ is
rather `regular' since one can compute at each point $x_0$ the same regularity exponent. This remark was actually the starting point of this work.\\

Indeed obviously the graph of the function has a very irregular behavior, and it has also some selfsimilarity properties. What can we say on the domain $\Omega=\{X=(x,y):y\leq F(x)\}$ under the graph of $F$ ?\\ 

{\it Denote in the following by $\ind_\Ome$ the characteristic function of $\Ome$, which takes the value 1 on $\Ome$ and 0 outside
$\Ome$.} \\

A first reflex is to compute fractal dimensions of the boundary $\partial\Omega$. The box
dimension of the graph can be derived by standard methods  (see Tricot \cite{Tricot}) and is exactly $\dim_B(\partial\Omega)=2-\alpha$.
Let us mention that Ciesielski \cite{Cie1, Cie2} proved results of this type for Schauder and Haar bases expansions in the case of more general families of functions.
Jaffard \cite{Ja53}, Kamont and Wolnik \cite{Kam67} obtained then general formulas  that allow to derive the box dimensions of the graphs
of arbitrary functions from their wavelet expansions.

For what concerns the Hausdorff dimension of the graph of $F$, as far as we know, the question is not solved yet in its all generality. It was proved by Ledrappier \cite{Ledrappier} in 1992 to be $2-\alpha$ in the special case where $a=2^{\alpha-1}$
is an Erd\"os number. By the results of Solomyak \cite{Solom} on Erd\"os numbers this amounts to have the computation for almost every $\alpha$ in [0,1].\\

Beside the computation of the box and Hausdorff dimension, which provide global quantities to describe the graph of the function, several methods were recently developed to classify fractal
 boundaries with the help of pointwise exponents. The idea was to be able to give a finer description of the geometry of the boundary, since the pointwise behavior was studied.
 In \cite{JaffMel}, Jaffard and M\'elot focused on the
 computation of the dimension of the set of points where $\ind_\Omega$  has a given $p$-exponent in the sense of
 Definition \ref{defi:pexponent}. In \cite{JaffHeur}, Jaffard and Heurteaux studied
 pointwise exponents more related to the geometry. These are the exponents we are actually interested in.\\

 Indeed denote by $meas$ the Lebesgue measure in $\R^d$ and $B(X,r)$ the $d$ dimensional open ball of center $X$ and radius $r>0$.
Jaffard and Heurteaux \cite{JaffHeur} gave the following definitions.

\begin{defi}\label{def:weakac}
Let $\Omega$ be a domain of $\R^d$ and let $X_0\in\partial\Omega$.
The point $X_0$ is weak $\alpha$-accessible in $\Omega$ if there exists $C>0$ and $r_0>0$ such that
\begin{equation}\label{eq:weakac}
\forall r\leq r_0 \qquad
meas(\Omega\bigcap B(X_0,r))\leq Cr^{\alpha+d}\;.
\end{equation}
The supremum of all the values of $\alpha$ such that (\ref{eq:weakac}) holds
 is called the weak accessibility exponent in $\Ome$ at $X_0$. We denote it by $E_\Omega^w(X_0)$.\\
\end{defi}

{\bf Example:}  Let
$0<\beta<1$ and $\Omega=\{X=(x,y)\in\R^2: |y|\leq |x|^\beta\}$. Denote $\Omega^c$
the complement of $\Omega$. Then one can easily check
that at each point $X_1\neq (0,0)$ of the boundary $\partial\Omega$ we have $E^w_\Omega(X_1)=0=E^w_{\Omega^c}(X_1)$
and at $X_0=(0,0) $ we have  $E^w_{\Omega^c}(X_0)=\frac{1}{\beta}-1 $ and $E^w_{\Omega}(X_0)=0$.
%

\begin{defi}\label{def:strongac}
Let $\Omega$ be a domain of $\R^d$ and let $X_0\in\partial\Omega$.
The point $X_0$ is strong $\alpha$-accessible in $\Ome$ if there exists $C>0$ and $r_0>0$ such that
\begin{equation}\label{eq:strongac}
\forall r\leq r_0 \qquad meas(\Omega\bigcap B(X_0,r))\geq Cr^{\alpha+d}\;.
\end{equation}
The infimum of all the values of $\alpha$ such that (\ref{eq:strongac}) holds is
called the strong accessibility exponent in $\Ome$ at $X_0$. We denote it by $E_\Omega^s(X_0)$.\\
\end{defi}
 The following proposition is given in \cite{JaffHeur}.

\begin{prop}\label{prop:definfsup}
Let $X_0 \in\partial\Omega$. Then

\begin{equation}\label{eq:liminf}
\begin{split}
d + E^w_\Omega(X_0) & = \liminf_{r \rightarrow 0}\frac{\log\left(meas(\Omega\cap B(X_0,r))\right)}{\log r}
\;, \\
d + E^s_\Omega(X_0) & = \limsup_{r \rightarrow 0}\frac{\log\left(meas(\Omega\cap B(X_0,r))\right)}{\log r}
\;.
\end{split}
\end{equation}
\end{prop}
Obviously $E^s_\Omega(X_0)\geq E^w_\Omega(X_0)$. We will see that thanks to our result one can prove that these two exponents can be different.\\

C.Tricot \cite{Tricot2}  proved that these exponents are related to local dimension computation. Let us mention, without entering too much the details, the relationship of this work \cite{Tricot2} with these exponents. Indeed the author focus on the formula 

\begin{equation}
H^\phi(E)=\liminf\limits_{\varepsilon\rightarrow 0}\{\sum\limits_{i\geq 0} \phi(E_i):E\subset \bigcup\limits_{i\geq 0}E_i, diam(E_i)\leq \varepsilon\}
\end{equation}
with $\phi:\mathbb{B}_E\rightarrow(0,\infty)$ some "set function" and $\mathbb{B}_E$ the set of closed balls centered on $E$.

Given an open set $V$ such that $E\subset \partial V$ the special choice of $$\phi_\alpha(B)=\frac{Vol(B\cap V)}{Vol(B)}diam(B)^\alpha$$
 lead to definitions of Hausdorff, exterior and interior dimensions, Packing, exterior and interior dimensions. \\
 
The following characterization, written for the setting we are interested in, holds

 \begin{thm}\cite{Tricot2}
 Let $\Omega$ be a bounded open set in $\R^d$ with boundary $\partial \Omega$ such that $meas(\partial\Omega)=0 $.
 Let $X_0\in \Omega$. Let $r>0$ and $$\alpha(B(X_0,r))=d-\frac{\log(meas(\Omega\bigcap B(X_0,r) ))}{\log(r)} \;.$$

 Then
 \begin{eqnarray*}
 \liminf\limits_{r\rightarrow 0}\alpha(B(X_0,r))&= dim_{int}(\{X_0\}) \qquad \mbox{and} \qquad
 \limsup\limits_{r\rightarrow 0}\alpha(B(X_0,r))&= Dim_{int}(\{X_0\})
 \end{eqnarray*}
 with $dim_{int} $ the Hausdorff interior dimension and $Dim_{int}$ the Packing Hausdorff dimension.

 \end{thm}

 We clearly have $dim_{int}(\{X_0\})=-E^w_\Omega(X_0)$ and $Dim_{int}(\{X_0\})=-E^s_\Omega(X_0)$.
  Let us stress that in the setting of Tricot $dim_{ext}(\{X_0\})=-E^w_{\Omega^c}(X_0)$ and
  $Dim_{ext}(\{X_0\})=-E^s_{\Omega^c}(X_0)$ with $\Omega^c$ the complementary of $\Omega$ in $\R^d$. We rather refer to \cite{Tricot2} for more details on local dimensions in their all generality.\\

 We will compute these quantities at the points of the boundary $\partial \Omega$ of
$\Omega=\{(x,y)\in\R^2:0\leq x\leq 1, 0\leq y\leq F(x)\}$, where  $F$ is the function defined by (\ref{expansion0}).
For that we will use the characterization of the maxima and minima done in \cite{Dub2}. This will yield the $p$-exponent at each point
 of $\ind_\Omega$. We will actually derive the fact that the set of local extrema of the function is fully characterized by the set of points where this $p$-exponent has a given value.\\
 We will also prove that the weak and strong accessibility exponents in $\Omega$ and $\Omega^c$ change from point to point
 on the graph $\partial \Omega$ of $F$. They also help to provide exact characterization of the sets of local maxima and local minima. Finally we will prove that there is a set of non trivial Hausdorff dimension such that the strong accessibility exponents in $\Omega$ and $\Omega^c$ are the same and the weak and strong accessibility exponents different. \\

Let us emphasize that this is to our knowledge the first time that the computation of these exponents was done in a nearly exhaustive study on a given example. The characterization we get for the set of extremas raise several questions: is it a general property ? Do other functions share it ? Could it lead to a finer classification of functions in H\"older classes ?  We would like to adress them in future works.\\

Let us come back now to our work.  The outline of the paper is the following.
In Section \ref{sec:defitheo} we set our main result.
 In Section \ref{sec:notations} some notations, preliminary remarks and technical lemmas, help us to prepare Section \ref{sec:proofs} where are the main proofs.

\vspace{24pt}

\section{Main results}\label{sec:defitheo}

\subsection{Statement of our main result}
Our goal is to prove the following Theorem.

\begin{mainthm}\label{th:results}
Let $a=2^{-\alpha}$ with $0<\alpha<1$ and $b=2$. Let $F$ be the function defined by (\ref{expansion0}).\\
Let $\Omega=\{X_0=(x,y)\in\R^2:0\leq x\leq 1, 0\leq y\leq F(x)\}$ and let $f=\ind_\Omega$.\\

Then at each point $X_0$ of $\partial\Omega$, the graph of $F$, we have
\begin{enumerate}
\item $u_f^p(X_0)=\frac{1}{p}\left(\frac{1}{\alpha}-1\right)$ if and only if $F(X_0)$ is a local extremum of $F$. Furthermore
\begin{enumerate}
\item $E^w_\Omega(X_0)=E^s_\Omega(X_0)=\frac{1}{\alpha}-1$ if and only if $F(X_0)$ is a local maximum of $F$. And in this case $E^w_{\Omega^c}(X_0)=E^s_{\Omega^c}(X_0)=0$.
\item $E^w_{\Omega^c}(X_0)=E^s_{\Omega^c}(X_0)=\frac{1}{\alpha}-1$ if and only if $F(X_0)$ is a local minimum of $F$. And in this case $E^w_\Omega(X_0)=E^s_\Omega(X_0)=0$.
\end{enumerate}
\item In the other cases where $F(X_0) $ is not a local extremum of $F$, we have $E^w_\Omega(X_0)=E^w_{\Omega^c}(X_0)=0 $.\\ \item Furthermore one can find a subset $D_\alpha\subset\partial\Omega$ such that for each $X_0\in D_\alpha$ $E^s_\Omega(X_0)=E^s_{\Omega^c}(X_0)=\frac{1}{\alpha}-1 $ and $E^w_\Omega(X_0)=E^w_{\Omega^c}(X_0)=0 $. The orthogonal projection of $D_\alpha$ on $[0,1]$ has the Hausdorff dimension $\alpha$.
\end{enumerate}
\end{mainthm}

\section{Useful notations and  results}\label{sec:notations}

\subsection{Lemmas for practical computation of the exponents}

From the computation of the weak accessibility exponent in $\Omega$ and $\Omega^c$ it is easy to derive the $p$-exponent. In \cite{JaffMel}, Jaffard
and M\'elot proved that $\ind_\Ome \in T_{\alpha/p}^p(X_0)$ if and only if either
$X_0$ is weak $\alpha$-accessible in $\Omega$ or $X_0$ is weak $\alpha$-accessible in $\Omega^c$. As a consequence we have
\BE \label{upgeqe}
p \;\; u^p_{\ind_\Ome } (X_0)  = \max(E^w_\Omega(X_0),E^w_{\Omega^c}(X_0))\;.
\EE

 We will also need the following lemma.

\begin{lem}\label{lem:Calpha}
Let $f : \R \rightarrow \R$ be in $C^\alpha(x_0)$ with $0< \alpha < 1$ and $\Omega$ be the domain
below (resp. above) the graph of $f$.
Consider $X_0=(x_0,f(x_0))$.
Then $X_0$ is strong $\frac1\al - 1$ accessible in both $\Ome$ and $\Ome^c$.
\end{lem}

{\bf Proof.} Suppose that $\Omega$ is the domain below the graph of $f$.
 Without any loss of generality, we can assume
that $X_0=(0,0)$. Let $r>0$.
Since $f$ is in $C^\alpha(0)$ and $0<\alpha<1$ then there exists a constant $C\geq 0$ such that in neighborhood of $0$
\begin{equation}
|f(x)| \leq C|x|^\alpha\;.
\end{equation}
Thus
\begin{equation}
-C|x|^\alpha \leq f(x) \leq C|x|^\alpha\;.
\end{equation}
Obviously $meas(\Omega^c\bigcap B(X_0,r))$ (resp.  $meas(\Omega \bigcap B(X_0,r))$) is greater than the
area $\mathcal{A} = C' \ds \int_0^r y^{1/\al} dy  = C" r^{1+1/\al}$
above (resp. below) the graph of $x\mapsto C|x|^\alpha $ and below (resp. above)
 the square of side $r$ and center $X_0$.\\
 The same results hold if $\Omega$ is the domain above the graph of $f$ (we have just to replace $\Ome$ by
 $\Omega^c$).


One of our goals for the points which are not extrema of $F$ will be to
find sequences of local maxima or minima such that
the following key-lemma proved in \cite{BS-H} holds.

\begin{lem}\label{lem:Calphaheur1}
Let $f : \R \rightarrow \R$ be in $C^\alpha(\R)$ and $\Omega$ be the domain below the graph of $f$.
Consider $X_0=(x_0,f(x_0))$. Suppose that there exist $c_\alpha>0$, a sequence $r_n$ of positive numbers, such that
 $r_n\rightarrow 0$ as $n\rightarrow+\infty$, $x_n\in ]x_0-r_n,x_0+r_n[$, and $n_0\in\N$, such that

\BE
\forall n\geq n_0\mbox{ }\mbox{ } f(x_n)=f(x_0)-c_\alpha r_n^\alpha\;.
\EE
Then $E^w_{\Omega^c}(X_0)=0$.
\end{lem}

{\bf Proof.} We can suppose that $x_n \in ]x_0-r_n,x_0[$ (the case $x_n \in ]x_0,x_0+r_n[$ is similar).
 Then by the mean value theorem we can find $b_n\in ]x_n,x_0[$ such that $f(b_n)=f(x_0)-r_n$.
Let $b'_n=\inf\{b_n\in ]x_n,x_0[ \;;\; f(b_n)=f(x_0)-r_n\}$.
Since $f$ is continuous we get $f(b'_n)=f(x_0)-r_n$ and $b'_n\in ]x_n,x_0[$.
It follows from the  definition of $b'_n$ and the mean value theorem that
$$ \forall t\in ]x_n,b'_n[ \quad f(t)<f(x_0)-r_n\;.$$
Thus $$]x_n,b'_n[\times[f(x_0)-r_n,f(x_0)+r_n]\subset B(X_0,r_n)\bigcap\Omega^c \;.$$
Therefore
\BE\label{eq:vol1}
meas(B(X_0,r_n)\bigcap\Omega^c )\geq 2|b'_n-x_n|r_n\;.
\EE
Since
\begin{eqnarray*}
f(x_n)-f(b'_n)&=&f(x_0)-c_\alpha r_n^\alpha-(f(x_0)-r_n)\\
&=&r_n-c_\alpha r_n^\alpha\\
&\simeq &-r_n^\alpha\;,
\end{eqnarray*}
in the sense that there exists a constant $C \geq 1$ such that for every $n$ we have $\ds \frac{1}{C}
r_n^\alpha \leq c_\alpha r_n^\alpha - r_n \leq Cr_n^\alpha$.\\
Since $f$ belongs to $C^\alpha(\R)$ we get
\BE
C|b'_n-x_n|^\alpha\geq |f(b'_n)-f(x_n)|\simeq r_n^\alpha\;.
\EE
Thus
\BE\label{eq:length}
|b'_n-x_n|\geq C' r_n\;.
\EE
Following (\ref{eq:vol1}) and (\ref{eq:length}) we get
\BE\label{eq:vol2}
meas(B(X_0,r_n)\bigcap\Omega^c )\geq Cr_n^2\;.
\EE
Since $meas(B(X_0,r)\bigcap\Omega^c )\leq meas(B(X_0,r))\leq r^2 $ for all $r\geq 0$
 we get
$$\lim\limits_{n\rightarrow+\infty}\frac{\log(meas(B(X_0,r_n)\bigcap\Omega^c ))}{\log(r_n)}=2\;.$$
Thus thanks to Proposition \ref{prop:definfsup} we have $E^w_{\Omega^c}(X_0)\leq 0$
which yields  $E^w_{\Omega^c}(X_0)=0$.\\

By replacing $f$ by $-f$ we also have the following result.

\begin{lem}\label{lem:Calphaheur2}
Let $f: \R \rightarrow \R$ be in $C^\alpha(\R)$ and $\Omega$ be the domain below the graph of $f$.
Consider $X_0=(x_0,f(x_0))$. Suppose that there exist $c_\alpha>0$,
  $r_n\rightarrow 0$ as $n\rightarrow+\infty$, $x_n\in ]x_0-r_n,x_0+r_n[$, and $n_0\in\N$, such that
\BE
\forall n\geq n_0\mbox{ }\mbox{ } f(x_n)=f(x_0)+c_\alpha r_n^\alpha\;.
\EE
Then $E^w_{\Omega}(X_0)=0$.
\end{lem}


\subsection{Dyadic expansions and approximation by dyadics}
We give some properties of the approximation of a point by the dyadics. Such
properties will be used later.\\

Let $x \in [0,1]$. Set $i_1(x), \cdots ,i_j(x),\cdots$  the binary digits of $x$,
i.e.
\begin{equation}\label{eq:dyadicexpan}
x=\ds \sum_{l=1}^{\infty} \frac{i_{l}(x)}{2^{l}}
\end{equation}

\begin{itemize}

\item Note that dyadic points, i.e points
$x=2^{-N} K$ with $K \in 2 \N +1$ are characterized by the fact that one can find $N>0$ such that
$i_N(x)=1$ and $i_n(x)=0$ for $n>N$, or equivalently $i_N(x)=0$ and $i_n(x)=1$ for $n>N$.\\

Furthermore for $n>N$ the number $x-2^{-n}$ is  dyadic. Since $2^{-n} = \ds \sum_{j=n+1}^\infty 2^{-j}$, then
$x - 2^{-n} = (\ds \sum_{j=1}^{N-1}i_j 2^{-j}) + 2^{-(N+1)}+ 2^{-(N+2)}+\cdots
+2^{-n}$.\\

On the other hand $x+2^{-n}$ has the simple expansion $\sum\limits_{j=1}^N\frac{i_j(x)}{2^j}+\frac{1}{2^n} $.\\

We will denote $\mathcal{D}$ the set of all dyadic points in $[0,1]$.\\

\item Let us come back to the general case with $x$ any point in $[0,1]$.

For each $j \in \NN$, define $K_j (=K_{j}(x))$ by

\BE\label{eq:approxdyad }
|K_{j}2^{-j} -x|=\inf_{k \in \NN}|k2^{-j} -x|\;.
\EE
Set
\[
r_{j}(x)= \frac{ \log |K_{j}2^{-j} - x|}{ \log 2^{-j}}.
\]
Define the rate of approximation of $x$ by dyadics as
\[
r(x)= \limsup_{j \mapsto \infty} r_{j}(x)\;.
\]
Since $|K_{j}2^{-j} - x| \leq 2^{-j}$, then for every $x$, we have $r(x) \geq
1$. If  $x$ is dyadic
then $r(x)=\infty$ (by taking the convention $\log 0 = - \infty$).
If  $x$ is normal (i.e. the frequency of ones (or zeros) in the binary
expansion of $x$ is equal to $1/2$) then $r(x)=1$.\\

\item If  $r(x)>1$, following the definition of $r(x)$, then for any $\delta>0$ such that $r(x)-\delta>1 $ one can find a subsequence $J_n\rightarrow +\infty $ for $n\rightarrow+\infty$ such that
\begin{equation}
r_{J_n(x)}\leq 2^{-J_n(r(x)-\delta)}\;.
\end{equation}
Let $J'_n=[J_n(r(x)-\delta)]$. We have then

\begin{equation}
K_{J_n}2^{-J_n}-2^{-J'_n}\leq x\leq K_{J_n}2^{-J_n}+2^{-J'_n}
\end{equation}
Thus, either $x$ belongs to the dyadic interval $[K_{J_n}2^{-J_n}-2^{-J'_n},K_{J_n}2^{-J_n}]$, and in this case it satisfies $i_{J_n+1}(x)=...=i_{J'_n-1}(x)=1 $, or it belongs to the other interval $[K_{J_n}2^{-J_n},K_{J_n}2^{-J_n}+2^{-J'_n}]$ and in this case it satisfies $i_{J_n+1}(x)=...=i_{J'_n-1}(x)=0 $.

In both cases let us notice that the binary expansion of $x$ contains chains of $0$ or $1$   whose length $J'_n-J_n\sim J'_n $ increases when $n\rightarrow+\infty$.

\end{itemize}
\vspace{24pt}

\subsection{Approximation by sequences of maxima of $F$}
We will see in the following that points in $[0,1]$ of the set
\BE
\mathcal{S}=\left\{k\in\N, N_0\in\N, \ds \frac{k}{2^{N_0}} + \frac{1}{3} \frac{1}{2^{N_0}},\ds \frac{k}{2^{N_0}} +
 \frac{2}{3} \frac{1}{2^{N_0}}\right\}
\EE
  will play a big role in this work, since they actually are the locations of the local maxima of the function $F$ (see below). Remark that they are characterized by the fact that for each $x\in \mathcal{S}$, one can find $j_0\in\N $ such that for $j\geq j_0$ we have $i_j(x)+i_{j+1}(x)=1 $.\\

As in the case of dyadic approximation we can define a rate of approximation by this kind of points.\\

Indeed let for $x\in [0,1]$
\BE\label{eq:approxmax }
|m_j -x|=\inf_{k \in \NN}\left\{\left|\frac{k}{2^{j}}+ \frac13 \frac{1}{2^{j}} -x\right|,\left|\frac{k}{2^{j}}+ \frac23 \frac{1}{2^{j}} -x\right|\right\}\;.
\EE
Define
\[
s_{j}(x)= \frac{ \log ( |m_j - x|}{ \log 2^{-j}}.
\]
Then the rate of approximation of $x$ by elements of $\mathcal{S}$  is given by
\[
s(x)= \limsup_{j \mapsto \infty} s_{j}(x)\;.
\]
Since $|m_j - x| \leq 2^{-j}$, then for every $x$, we have $s(x) \geq
1$.\\

In the case of dyadic numbers, we have $s(x)=1$. But
remark that in other non trivial cases there is no obvious relationship between $s(x)$ and $r(x)$.
Indeed one can check on the following examples that $s$ and $r$ can take independently any value.

\begin{itemize}
\item Let $x=\sum\limits_{j=1}^{+\infty}\frac{i_j(x)}{2^j}$ with $i_{3k+1}(x)=0=i_{3k+2}(x) $ and $i_{3k+2}(x) =1 $ for all $k\in\N$. Then we have $r(x)=s(x)=1 $.
\item Let $u>1$. Then $x=\sum\limits_{n=0}^{+\infty}2^{-[u^n]}$ with $[u^n] $ the integer part of $u^n$. We have $r(x)=u$ whereas $s(x)=1 $.
\item Let $u>1$. Then $x=\sum\limits_{n=1}^{+\infty}2^{-2n}-\sum\limits_{n=0}^{+\infty}2^{-2[u^n]}$. We have $s(x)=u$ whereas $r(x)=1 $.
\item Let $u>1$ and $s>1$. Let $x=\sum\limits_{n=1}^{+\infty}2^{-2[s^{n-1}u^n]}+\sum\limits_{k=[s^{n-1}u^n]+1}^{[s^{n}u^{n+1}]-1}2^{-2k}$. Then $r(x)=u$ and $s(x)=s$.

\end{itemize}

\vspace{24pt}

\subsection{The shift operator}

Since $0 < \al < 1$  it is easy to check that we obtain from (\ref{expansion}) with $a=2^{-\alpha}$ and $b=2$
\BE \label{selfFeg}
  F(x) = \sum_{j=0}^{\infty} \;
\sum_{
      i = (i_1, \cdots, i_j) \in \{0,1\}^j
      } 2^{-\al j} \;
 \La \left
  (2^{j} x - 2^{j-1}i_1- \cdots - 2 i_{j-1}-i_j
\right).
\EE
The term of (\ref{selfFeg}) corresponding to $j=0$ is $\La(x)$.
But, the function $\La$ is supported in $[0,1]$, therefore $F$ vanishes outside
$[0,1]$ and for $x \in [0,1]$
\BE \label{ch5eq:1.2}
  F(x) = \sum_{j=0}^{\infty}
 2^{-\al j}  \,\,
\La \left
  (2^{j} x - 2^{j-1}i_1(x)- \cdots - 2 i_{j-1}(x)-i_j(x)
\right)\;.
\EE
For dyadic rationals $x$, $x=2^{-N} K$ with $K \in 2 \NN + 1$, as we already said it, there exist
two binary expansions, one such that $i_N(x)=1$ and $i_n(x)=0$ for $n>N$,
and another one such that $i_N(x)=0$ and $i_n(x)=1$ for $n>N$.
The two right-hand sides of (\ref{ch5eq:1.2}) corresponding to the two
choices of $i(x)$ give identical results.

Denote by $\tau$ the shift operator
\[
\tau x = \sum_{l=2}^{\infty}i_{l}(x)2^{-l+1} =
\sum_{l=1}^{\infty}i_{l+1}(x)2^{-l}\;.
\]
Observe that
$$
\tau x =  \left \{
\begin{tabular}{ll}
$2x$ & if $x \in [0,1/2[$\\
$2x-1$ & if $ x \in [1/2,1]$.
\end{tabular}
\right.
$$
Hence
\begin{equation}\label{eq:Falpha}
F(x) =  \sum_{j=0}^{\infty}  2^{-\al j}\,\, \La\left(\tau^{j}x\right)
\end{equation}
and
\[
\tau^j x = \sum_{l=1}^{\infty}i_{l+j}(x)2^{-l}\;.
\]

Our selfsimilar function  is of the form
$ F(x) = \ds \sum_{j=0}^{\infty} \ds \sum_{k=0}^{2^j - 1} C_{j,k} \La(2^{j}x-k) $
with
$$\begin{array}{lll}
C_{j,k} &=& 2^{-\al j}\mbox{ if }j\neq 0,k\neq 0\\
C_{0,0}&=&1
\end{array}
$$

\vspace{24pt}

For $n \geq 1$,  denote
\begin{equation}\label{eq:Fnalpha}
F_n(x) =  \sum_{j=0}^{n}  2^{-\al j}\,\, \La\left(\tau^{j}x\right)\;.
\end{equation}
Remark that $F_{n}$ is affine on intervals of type $I_{n+1}=\left[\frac{k}{2^{n+1}},\frac{k+1}{2^{n+1}}\right] $.
Remark also that if $t \in [0,1]$, then $\La(t) =  (-1)^{i_1(t)}t + i_1(t)$. So, if
$t' \in [0,1]$ and $i_1(t) = i_1(t')$, then
    $\La(t)-\La(t')= (-1)^{i_1(t)}(t-t')$.\\
It follows  that if $\ds \frac{k}{2^{n+1}} = \sum_{j=1}^{n+1} \frac{i_j}{2^{j}}$ then
 the slope of $F_{n}$ at any point $x$ of the interval $]\frac{k}{2^{n+1}},\frac{k+1}{2^{n+1}}[ $
  is exactly
\begin{equation}\label{eq:slope1}
C_{n} = C_{n}(x) = \sum\limits_{j=0}^{n}(-1)^{i_{j+1}(x)}2^{(1-\alpha)j}=
\sum\limits_{j=0}^{n}(-1)^{i_{j+1}}2^{(1-\alpha)j}\;.
\end{equation}

\subsection{Extrema of $F$}
We will need the following characterization of the extremas of $F$ proved in \cite{Dub1} and  \cite{Dub2}. Let us start with the local and global minima.\\

\begin{prop}\label{prop:minmax}
Let $0<\alpha<1$ and $F$ the function defined by (\ref{eq:Falpha}), then
\begin{itemize}
\item 0 and 1 are the abscissas of the global minima of $F$.
\item The dyadic points are the abscissas of the minima of $F$ and furthermore
\begin{equation}
\min\limits_{x\in I_{N}}F(x)=\min\left[F_{N-1}\left(\frac{k}{2^N}\right),F_{N-1}\left(\frac{k+1}{2^N}\right)\right]
\end{equation}
with $I_{N}=\left[\frac{k}{2^N}, \frac{k+1}{2^N}\right]$.
\end{itemize}
\end{prop}
In the case of the maxima, the statement of the result is slightly more technical. We need the following proposition of \cite{Dub1} using the same notations as previously.

\begin{prop}\label{prop:minmaxM}
Let $0<\alpha<1$ and $F$ the function defined by (\ref{eq:Falpha}). Let $t=2^{1-\alpha}$ and $X(p)$ the list of positions where $F(x) + px$ attains its
maximum on $[0, 1]$. Let $\mathcal{M}(p)$ be the maximum on $[0,1]$ of $F(x)+px $. Then
\begin{itemize}
\item $X\left(\frac{-(t^N-1)}{t-1}\right)=\left\{\frac{1}{3} \frac{1}{2^N},\frac{2}{3}
\frac{1}{2^N}\right\}$ for $N=0,1,...$
\item $X(p)=\left\{\frac{1}{3} \frac{1}{2^N}\right\}$ if $\frac{-(t^{N+1}-1)}{t-1}< p<  \frac{-(t^{N}-1)}{t-1}$.
\item $X(p)=1-X(-p)$ for all $p\neq 0$.

%
%
\item $\max\limits_{x\in I_N} F(x)=F_{N-1}\left(\frac{k}{2^N}\right)+
2^{-N\alpha}\mathcal{M}\left(\frac{C_{N-1}\left(\frac{2k+1}{2^{N+1}}\right)}{t^{N}}\right)$.
\end{itemize}
\end{prop}

The following proposition is a consequence of the previous one.

\begin{prop}\label{prop:minmax2}
Let $0<\alpha<1$ and $F$ the function defined by (\ref{eq:Falpha}). Then
\begin{itemize}
\item $1/3$ and $2/3$ are the abscissas of the global maxima of $F$.
\item The abscissas of the local maxima of $F$ are the points of $\mathcal{S}$.
\end{itemize}
\end{prop}

\subsection{Approximation of slopes of $F_n$}\label{sec:slopes}
Suppose first we have some informations about the dyadic expansion of $x$.
Then we have the following Lemma, which helps to control the behavior of the slopes of the affine function $F_{n-1}$.
\begin{lem}\label{prop:specialslopes}

\begin{enumerate}
\item\label{it:dyadslope} Let $x=\frac{K}{2^N}$ be a dyadic number. Then one can find $N_0>N$, $A>0$ and $B>0$
depending only on $x$ such that if $n\geq N_0 $ then
\begin{equation}\label{eq:dyadslope1}
\forall \; y\in] x,x+2^{-n}[   \quad A 2^{(1-\alpha) n}\leq C_{n-1}(y)\leq B2^{(1-\alpha) n}
\end{equation}
 and
\begin{equation}\label{eq:dyadslope2}
\forall \; y\in] x-2^{-n},x[ \quad   -A 2^{(1-\alpha) n}\geq C_{n-1}(y)\geq -B2^{(1-\alpha) n}\;.
\end{equation}

\item\label{it:maxslope} Let $x$ be the abscissa of a local maximum of $F$. Then one can find $J_0>0$, $A>0 $ and $B>0$ such that for $n\geq J_0$
\begin{equation}\label{eq:maxslope}
\begin{split}
i_{n-1}(x)+i_n(x)&=1\;,\\
C_{n-1}(x)C_n(x)&<0\\
A 2^{(1-\alpha) n}\leq (-1)^{i_{n+1}(x)}C_{n-1}(x)&\leq B2^{(1-\alpha) n}\;.
\end{split}
\end{equation}

\item\label{it:r1slope} Let $x$ be a non dyadic point such that $r(x)>1$. Then one can find two subsequences $J_n$ and $J'_n$ with $
\frac{J'_n}{J_n}>1$ for all $n$, such that $i_{J_n}(x)=i_{J'_n+1}(x)$ and $i_j(x)+i_{J_n}(x)=1$ for $J_n<j<J'_n+1$.
Furthermore one can find $J_0>0$, $A>0$ and $B>0$ such that for $n>J_0$
\begin{equation}\label{eq:r1slope}
A 2^{(1-\alpha) J'_n}\leq (-1)^{i_{J_n+1}(x)}C_{J'_n-1}(x)\leq B2^{(1-\alpha) J'_n}\;.
\end{equation}

\item\label{it:s1slope} Let $x$ be a non dyadic point such that $s(x)>1$. Then one can find two subsequences $J_n$ and $J'_n$ with $
\frac{J'_n}{J_n}>1$ for all $n$, such that $i_j(x)+i_{j+1}(x)=1$ for $J_n<j<J'_n$ and $i_{J'_n}(x)=i_{J'_n+1}(x) $. Furthermore one can find $J_0>0$, $A>0$ and $B>0$ such that for $n>J_0$
\begin{equation}\label{eq:s1slope}
A 2^{(1-\alpha) J'_n}\leq (-1)^{i_{J'_n}(x)}C_{J'_n-1}(x)\leq B2^{(1-\alpha) J'_n}\;.
\end{equation}

\end{enumerate}
\end{lem}

{\bf Proof.}
\begin{itemize}

\item Case \ref{it:dyadslope}: the idea is very simple since it is a direct computation.

Indeed following (\ref{eq:slope1}) we have for $y\in ]x,x+2^{-n}[$
\begin{eqnarray*}
 C_{n-1}(y)& = &
  \sum_{j=0}^{n-1} (-1)^{i_{j+1}(y)} 2^{-\al j} 2^j \\
&=& \sum_{j=0}^{N-2} (-1)^{i_{j+1}(y)} 2^{j(1-\al)}
-  2^{(N-1)(1-\al)} +  \sum_{j=N}^{n-1}2^{j(1-\al)}\\
&=& \sum_{j=0}^{N-2} (-1)^{i_{j+1}(y)} 2^{j(1-\al)}
+2^{N(1-\al)}(1-2^{\alpha-1})+2^{(N+1)(1-\alpha)}\frac{2^{(n-N-1)(1-\alpha)}-1}{2^{(1-\alpha)}-1}\;.
\end{eqnarray*}
The second equation with $y\in]x-2^{-n},x[$ can be computed in the same way, up to a change of signs.\\

Thus one can find $N_0>N$, $A>0$  and $B>0$ such that (\ref{eq:dyadslope2}) holds for $n>N_0$.

\item  Case \ref{it:maxslope}: it is enough to remark that $x$ has the following binary expansion
\begin{equation}\label{eq:binexmax}
x=\frac{k_{N_0}}{2^{N_0}}+\sum\limits_{l=0}^\infty\frac{1}{2^{2l+1+N_0}}\;.
\end{equation}
As a consequence of Proposition \ref{prop:minmax2}, the same kind of computation yields Case 2.

\item Case \ref{it:r1slope}: since $r(x)>1$, for any $\delta>0$ one can find two subsequences $J_n$ and $J'_n$
such that $i_{J_n+1}(x)=...=i_{J'_n-1}(x)$ and $|x-K_{J_n}2^{-J_n}|\leq 2^{-J'_n} $ with $J'_n=[(r(x)-\delta) J_n] $.\\
 Suppose first eventually up to a small change of definition of $J'_n$ that $i_{J_n}(x)=1=i_{J'_n+1}(x) $
 and $i_{J_n+1}(x)=...=i_{J'_n-1}(x)=i_{J'_n}(x)=0$. Then with the same kind of computation as in Case \ref{it:dyadslope} one gets

\begin{equation}\label{eq:r1slopesansabs}
A 2^{(1-\alpha) J'_n}\leq C_{J'_n-1}(x)\leq B2^{(1-\alpha) J'_n}\;.
\end{equation}

In the other case $i_{J_n+1}(x)=...=i_{J'_n-1}(x)=i_{J'_n}(x)=1$, the sign of the slope will be changed.

\item Case \ref{it:s1slope}: this follows exactly the same ideas than previously.
Since $s(x)>1$ for any $\delta>0$ one can find two subsequences $J_n$ and $J'_n$ such that
for all $J_n< j<J'_n$ $i_j(x)+i_{j+1}(x)=1$ and $i_{J_n}(x)=i_{J_{n}-1}(x)$, $i_{J'_n}(x)=i_{J'_n+1}(x) $.
Then with the same kind of computation as in Case \ref{it:dyadslope} one gets

\begin{equation}\label{eq:s1slopesansabs}
A 2^{(1-\alpha) J'_n}\leq (-1)^{i_{J'_n}(x)}C_{J'_n-1}(x)\leq B2^{(1-\alpha) J'_n}\;.
\end{equation}

Hence the result.
\end{itemize}
\vspace{10pt}
If we don't have any further information on $x$, the following Lemma will be useful.\\


\begin{lem}\label{prop:slopes}
 Let $x\in [0,1]$ be a non dyadic number.
 \begin{enumerate}
\item\label{it:slope1}\cite{Gau} Then there exists $\delta>0$ and $\delta'=\frac{1}{2^{1-\alpha}-1}$ such that for all $n\in\N$ one can find $J_n\geq n$ such that
\begin{equation}\label{eq:slopes}
\delta' 2^{J_n(1-\alpha)}>\left|C_{J_n-1}(x)\right|>\delta 2^{J_n(1-\alpha)}\;.
\end{equation}
\item\label{it:slope2} If $x\notin\mathcal{S}$ then there exists $\delta>0$ such that for all $n\geq 0$ there exists $J_n\geq n$ such that
(\ref{eq:slopes}) holds and
\begin{itemize}
\item either $\left[C_{J_n-1}(x)>0\mbox{ and }i_{J_{n}+1}(x)=0\right],$
\item or $\left[C_{J_n-1}(x)<0\mbox{ and }i_{J_{n}+1}(x)=1\right]$.
\end{itemize}
\end{enumerate}
\end{lem}

{\bf  Proof. }
\begin{enumerate}
\item The upper bound is a straightforward computation.\\

Suppose the contrary, i.e for all $\delta>0$ one can find $N\in\N$ such that for all $n\geq N$

\begin{equation}
-\delta 2^{(1-\alpha)n}\leq C_{n-1}(x)\leq \delta 2^{(1-\alpha)n}\;.
\end{equation}
If we suppose without loss of generality that $i_{n+1}(x)=0 $ then at step $n$

\begin{equation}
\begin{split}
-\delta 2^{(1-\alpha)n}\leq &C_{n-1}(x)\leq \delta 2^{(1-\alpha)n}\\
-\delta 2^{(1-\alpha)n}+2^{(1-\alpha)n}\leq &C_{n}(x)\leq \delta 2^{(1-\alpha)n}+2^{(1-\alpha)n}\\
-\frac{\delta}{2^{1-\alpha}}+\frac{1}{2^{1-\alpha}}\leq &\frac{C_{n}(x)}{2^{(1-\alpha)(n+1)}}\leq \frac{\delta}{2^{1-\alpha}}+\frac{1}{2^{1-\alpha}}\;.
\end{split}
\end{equation}
It is enough to choose $\delta$ such that $-\frac{\delta}{2^{1-\alpha}}+\frac{1}{2^{1-\alpha}}> \delta $ to have a contradiction.\\

\item Suppose the contrary, i.e there exists $x\notin\mathcal{S}$ and that for all $\beta>0$, there exists $N\in \N$, such that for all $n\geq N$
\begin{enumerate}
\item\label{it:beta} either $|C_{n-1}(x)|\leq\beta 2^{J_n(1-\alpha)} $,
\item\label{it:negslope} or $\left[C_{n-1}(x)< 0\mbox{ and } i_{n+1}(x)=0 \right]$,
\item\label{it:posslope} or $\left[C_{n-1}(x)> 0\mbox{ and } i_{n+1}(x)=1 \right]$.
\end{enumerate}

Remark first that the points of $\mathcal{S}$ satisfy exactly (\ref{it:negslope})
 and (\ref{it:posslope}). Indeed for $x\in\mathcal{S}$, and assuming that $i_{n_0+1}(x)=1 $, $x$ has a binary expansion (\ref{eq:binexmax}).\\

 Thus following (\ref{eq:slope1}) the slope $C_{n-1}(x) $ satisfies

 \begin{eqnarray*}
C_{n-1}(x)&=\sum\limits_{j=0}^{N_0-1}(-1)^{i_{j+1}(x)}2^{(1-\alpha)j}
+\sum\limits_{j=N_0}^{n-1} 2^{(1-\alpha)j}(-1)^{i_{j+1}(x)}\\
&=\sum\limits_{j=0}^{N_0-1}(-1)^{i_{j+1}(x)}2^{(1-\alpha)j}
+2^{N_0(1-\alpha)}
\frac{(-2^{(1-\alpha)})^{n-N_0}-1}{1+2^{1-\alpha}}\;.
\end{eqnarray*}
Hence, for $n$ large enough (\ref{it:negslope})
 and (\ref{it:posslope}) are satisfied.\\

 Our goal is thus to prove that if we choose $\beta$ small enough then only (\ref{it:negslope})
 and (\ref{it:posslope}) can be satisfied, which will lead to the fact that $x\in\mathcal{S}$, and thus to a contradiction.

Let start by the following special cases.
\begin{itemize}
\item We claim that if one can find $k$ large enough such that $C_{k-1}(x)=0 $ then $x\in\mathcal{S}$, which is a contradiction.\\
Let us prove this claim.\\

We will need the following sequence: let for $n\in\N^\star$
\begin{equation}\label{eq:dn}
d_n=2^{-(1-\alpha)}\sum\limits_{j=0}^n(-1)^j2^{-j(1-\alpha)}=\frac{2^{-(1-\alpha)} }{1+2^{-(1-\alpha)}}\left(1-(-1)^{n+1}2^{-(n+1)(1-\alpha)}\right)\;.
\end{equation}

We have clearly $d_n\geq d_1> 0$ for all $n\geq 1$.

Choose $\beta\leq \frac{d_1}{2}$ and $N$ such that the hypothesis are satisfied.\\

Suppose that $k\geq N+1$ is such that $C_{k-1}(x)=0 $. Then

 $|C_k(x)|=2^{k(1-\alpha)}=2^{-(1-\alpha)}2^{(k+1)(1-\alpha)} $.\\
Remark that $\beta<d_1\leq 2^{-(1-\alpha)} $ thus $|C_k(x)|>\beta 2^{(k+1)(1-\alpha)} $.\\

 Suppose without lost of generality that $C_k(x)>0$ (the case $C_k(x)<0$ is symetrical and can be proved in exactly the same way). Thus $i_{k+2}(x)=1 $ and
 \begin{equation}
 \begin{split}
 C_{k+1}(x)=-2^{(k+1)(1-\alpha)}+2^{k(1-\alpha)}&=-2^{(k+2)(1-\alpha)}\left(2^{-(1-\alpha)}-2^{-2(1-\alpha)}\right)\\
 &<-\beta 2^{(k+2)(1-\alpha)}\;.
 \end{split}
 \end{equation}


Let us prove by induction on $n$ that for all $n\in\N^\star $
$$|C_{k+n}(x)|= d_n2^{(k+n+1)(1-\alpha))},\mbox{}(-1)^{n}C_{k+n}(x)>0 \quad \mbox{ (P) }.$$

 We just proved that (P) is true for $n=1 $.\\

Suppose that for $n\in\N^\star$ (P) is true. Suppose without lost of generality that $C_{k+n}(x)>0 $ (the case
$C_{k+n}(x)<0$ is symetrical and can be proved in exactly the same way ). Thus $i_{k+n+2}(x)=1 $ and

  \begin{equation}
 \begin{split}
 C_{k+n+1}(x)&=C_{k+n}(x)-2^{(n+k+1)(1-\alpha)}=d_n2^{(k+n+1)(1-\alpha))}-2^{(n+k+1)(1-\alpha)}\\
&=-(-2^{-(1-\alpha)}d_n+2^{-(1-\alpha)})2^{(k+n+2)(1-\alpha)}\;.
 \end{split}
 \end{equation}

Since  $d_n$ satisfies exactly $d_{n+1}=-2^{-(1-\alpha)}d_n+2^{-(1-\alpha)}$, we have the result and (P) is satisfied at level $n+1$.\\

Thus for all $n\in\N^\star$ (P) is true. Remind that since  $\beta<d_1\leq d_n$ for all $n\in\N^\star$, this implies that for all $n\in\N^\star$ $i_{k+n+2}(x)+i_{k+n+3}(x)=1 $,
which is exactly the characterization of the points in $\mathcal{S}$, and is in contradiction with the hypothesis $x\notin\mathcal{S}$.\\

In the following we will always keep the hypothesis $0<\beta\leq \frac{d_1}{2} $ so that for $n$ large enough we have always $C_n(x)\neq 0 $.

\item We now consider the case where $|C_{n-1}(x)|$ is close to the value of $\beta 2^{n(1-\alpha)}$ and prove that this yields that $x\in\mathcal{S}$, and thus a contradiction.\\

Let $0<\beta\leq \frac{d_1}{2}$, and $\beta'>0$ whose value will be precised later on.\\

Suppose $n\geq N$ is such that $\beta' 2^{n(1-\alpha)} >C_{n-1}(x)>\beta 2^{n(1-\alpha)} $\\
(the case $C_{n-1}(x)<-\beta 2^{n(1-\alpha)}$ can be done exactly in the same way). Then $i_{n+1}(x)=1 $ and $C_{n}(x)=C_{n-1}(x)-2^{n(1-\alpha)}$, hence

\begin{equation}
\begin{split}
\beta 2^{n(1-\alpha)}-2^{n(1-\alpha)}\leq& C_n(x)\leq \beta'2^{n(1-\alpha)}-2^{n(1-\alpha)}\\
(\beta-1)2^{-(1-\alpha)}2^{(n+1)(1-\alpha)}&\leq C_n(x)\leq (\beta'-1)2^{-(1-\alpha)} 2^{(n+1)(1-\alpha)}\;.
\end{split}
\end{equation}
%
Choose $\beta'$ such that $(\beta'-1)2^{-(1-\alpha)}<-\beta$, hence $\beta<\beta'<\frac{1-\beta}{2^{-(1-\alpha)}} $,
 which is possible since $\beta\leq \frac{d_1}{2}<\frac{2^{-(1-\alpha)}}{1+2^{-(1-\alpha)}} $.\\

This yields $i_{n+2}(x)=0 $. Thus
\begin{equation}
\begin{split}
C_{n+1}(x)&=C_{n-1}(x)+(-2^{n(1-\alpha)}+2^{(n+1)(1-\alpha)})\\
&>(2^{-(1-\alpha)}-2^{-2(1-\alpha)})2^{(n+2)(1-\alpha)}>\beta 2^{(n+2)(1-\alpha)}\;.
\end{split}
\end{equation}

Let us prove by induction on $k$ that
for all $k\in\N$, $$(-1)^{k+1}C_{n+k}(x)>\beta 2^{(n+k+1)(1-\alpha)} \quad \mbox{ (Q) }\;.$$
We just prove that the case $k=0 $ is true.\\

Suppose one can find $k\in\N$ such that (Q) is true for all $0\leq k'\leq k$.\\
Let us prove that it is true at $k+1$. Without lost of generality suppose $C_{n+k}(x)>0 $, thus $i_{n+k+2}(x)=1 $.\\
We have
\begin{equation}
\begin{split}
C_{n+k+1}(x)&=C_{n+k-1}(x)+2^{(n+k)(1-\alpha)}-2^{(n+k+1)(1-\alpha)}\\
&<-(2^{-(1-\alpha)}-2^{-2(1-\alpha)})2^{(n+k+2)(1-\alpha)}<-\beta 2^{(n+k+2)(1-\alpha)}\;.
\end{split}
\end{equation}
This proves that (Q) is true at $k+1$.\\

 Thus by induction (Q) is true for all $k\in\N $. This means that for all $k\in\N $ $i_{n+k+2}(x)+i_{n+k+3}(x)=1$, and thus $x\in\mathcal{S}$. Hence the contradiction.\\

\item We now study the case where  $|C_{n-1}(x)|\leq \beta 2^{n(1-\alpha)} $ and prove that if we choose $\beta$ small enough then it will lead to $x\in\mathcal{S}$.\\

Indeed let $\beta\leq\inf(\frac{d_1}{2},\frac{d'_1}{2}) $ with
\begin{equation}\label{eq:d'1}
d'_1=\frac{2^{-(1-\alpha)}-2^{-2(1-\alpha)}}{1+2^{-2(1-\alpha)}}\;.
\end{equation}

And suppose $n\geq N$ such that $|C_{n-1}(x)|\leq \beta 2^{n(1-\alpha)} $. Suppose $i_{n+1}(x)=0 $ without lost of generality. Thus we have
\begin{eqnarray*}
&C_n(x)&=C_{n-1}(x)+2^{n(1-\alpha)}\\
(-\beta+1)2^{-(1-\alpha)}2^{(n+1)(1-\alpha)}&\leq C_n(x)&\leq  (\beta+1)2^{-(1-\alpha)}2^{(n+1)(1-\alpha)}\;.
\end{eqnarray*}

Remark that with the choice of $\beta$ we made, we have on one hand $\beta<(-\beta+1)2^{-(1-\alpha)}$ and
on the other hand $\beta<(\beta+1)2^{-(1-\alpha)}<\frac{1-\beta}{2^{-(1-\alpha)}}$. Thus following the previous
result using $\beta'=(\beta+1)2^{-(1-\alpha)}$, $x\in\mathcal{S}$ and we have a contradiction.

\item We consider the case where $C_{n-1}(x)>0$ and $C_n(x)<0$ for $n$ large enough under the previous range of values of $\beta$.\\

Let $\beta\leq\inf(\frac{d_1}{2},\frac{d'_1}{2}) $ (recall that $d_1$ is defined by (\ref{eq:dn})  and $d'_1$ by (\ref{eq:d'1})).\\
And suppose that for $n\geq N$ we have $C_{n-1}(x)>0$ and $C_n(x)<0$.\\
Following the previous case we have $C_{n-1}(x)>\beta 2^{n(1-\alpha)} $ and $C_n(x)<-\beta 2^{(n+1)(1-\alpha)}$.
 Thus $i_{n+1}(x)=1 $ and $i_{n+2}(x)=0$.\\

Then
\begin{equation}
\begin{split}
C_{n+1}(x)&=C_{n-1}(x)-2^{n(1-\alpha)}+2^{(n+1)(1-\alpha)}\\
&>\beta 2^{(n+2)(1-\alpha)}
\end{split}
\end{equation}

since by definition of $\beta $ and $d_1 $ we have $2^{-(1-\alpha)}-2^{-2(1-\alpha)}\geq d_1>\beta $. Thus $i_{n+3}(x)=1$.\\
We have $C_{n+2}(x)=C_{n}(x)+2^{(n+1)(1-\alpha)}-2^{(n+2)(1-\alpha}<-\beta 2^{(n+3)(1-\alpha)}$.\\

 A proof by induction exactly in the same way as previously yields that for $k\geq 0$ we have $C_{n+2k+1}(x)>0 $ and $C_{n+2k}(x)<0 $, thus $i_{n+2k+2}(x)+i_{n+2k+3}(x)=1 $ for all $k\in\N$ and we have $x\in\mathcal{S}$, hence a contradiction.\\
%
\end{itemize}
We will now go the main proof, taking into account what we just proved.\\

In the following we will consider $\delta>0$ and $J_n$ defined as in Point \ref{it:slope1},  $\beta=\inf(\frac{d_1}{2},\frac{\delta}{2},\frac{d'_1}{2})$ and $n$ such that $J_n\geq N$. Thus for all $n\in\N$ $\frac{|C_{J_n-1}(x)|}{2^{J_n(1-\alpha)}}>\delta>\beta$.\\

Suppose $C_{J_{n}-1}(x)>0 $. This means that $C_{J_n}(x)=C_{J_{n}-1}(x)-2^{J_n(1-\alpha)}<C_{J_{n-1}}(x) $. The only case we want to consider is $C_{J_n}(x)>\beta 2^{(J_n+1)(1-\alpha)}$ since for all the other cases the previous points yield $x\in\mathcal{S}$.\\

 Thus $i_{J_n+2}(x)=1 $. It is clear that one can find  $k\in\N$ such that for all $0\leq k'\leq k$ $C_{J_n+k'}(x)>\beta 2^{(J_n+k'+1)(1-\alpha)}$ and $i_{J_n+k'+2}(x)=1$ and $C_{J_n+k+1}(x)<\beta 2^{(J_n+k+2)(1-\alpha)}$.\\

Hence either $|C_{J_n+k+1}(x)|\leq  \beta 2^{(J_n+k+2)(1-\alpha)}$ and $x\in\mathcal{S}$, or $C_{J_n+k+1}(x)<-\beta 2^{(J_n+k+2)(1-\alpha)} $ and since $C_{J_n+k}(x)>\beta 2^{(J_n+k+1)(1-\alpha)} $ we have also the result.\\

In all cases we proved that Points (\ref{it:beta}), (\ref{it:negslope}), (\ref{it:posslope}) lead to $x\in\mathcal{S}$, which is a contradiction. Hence the Lemma.

\end{enumerate}

\vspace{24pt}

\section{Computation of weak and strong accessible exponents}\label{sec:proofs}
\subsection{Case of dyadic points}

\vspace{24pt}
We will prove the following proposition.

\begin{prop} \label{xdyadic}
If $x$ is a dyadic point, and $X=(x,F(x))$ then
\BE \label{Ewdyadic}
E^w_{\Ome^c}(X) = \frac{1}{\al} - 1\mbox{ and }E^w_{\Ome}(X) =0
\EE
\BE \label{Esdyadic}
E^s_{\Ome^c}(X) = \frac{1}{\al} - 1\mbox{ and }E^s_{\Ome}(X) =0
\EE
\BE \label{updyadic}
u^p_f(X) = \ds \frac{1}{p} (\ds \frac{1}{\al} - 1)\;.
\EE
\end{prop}

If $x$ is a dyadic point, i.e. $x=2^{-N} K$ with $K \in 2 \NN +
1$, we consider its  binary expansion  in which $i_N(x)=1$ and
$i_n(x)=0$ for $n>N$.
For $n>N$ the number $x-2^{-n}$ is  dyadic. Since $2^{-n} = \ds \sum_{j=n+1}^\infty 2^{-j}$ then
$x - 2^{-n} = (\ds \sum_{j=1}^{N-1}i_j 2^{-j}) + 2^{-(N+1)}+ 2^{-(N+2)}+\cdots
+2^{-n}$.
On the other hand $x+2^{-n}$ has the simple expansion
 $(\ds \sum\limits_{j=1}^N\frac{i_j(x)}{2^j})+\frac{1}{2^n} $.\\
Remark that $F_{N-1}(x)=F_{n-1}(x)=F(x)$ for $n> N$ and $F_{n-1}(x+2^{-n})=F(x+2^{-n}) $.\\

%

 Any point $y$ in the interval $]x,x+2^{-n}[$ satisfies the expansion $i_N(y)=i_{N+1}(y)=...=i_{n-1}(y)=i_n(y)=0$.\\
It follows that
\begin{eqnarray*}
F(x+2^{-n}) - F(x)= F_{n-1}(x+2^{-n}) - F_{n-1}(x) = 2^{-n} C_{n-1}(y)
  \end{eqnarray*}
 with $y$ any of the points of the interval $]x,x+2^{-n}[$.

 Following Lemma \ref{prop:specialslopes} and Case \ref{it:dyadslope} there exist two constants $A>0$ and $B>0$
 and $J_0 \geq N$ (which depend only on the given dyadic point
$x$) such that
\BE \label{F+J-Fsoit}
\forall n \geq J_0 \quad
 A 2^{-\al n}\leq F(x+2^{-n}) -  F(x)\leq B2^{-\alpha n}\;.
\EE
Thus we have $F(x+2^{-n})-F(x)\geq A.2^{-\alpha n}$.

%

On the other hand, following remarks of Section \ref{sec:notations3-2}, for any $y\in ]x-2^{-n},x[$ we have $i_N(y)=1=...=i_{n}(y) $. Thus
\begin{equation}
 F(x-2^{-n}) - F(x)  =  F_{n-1}(x-2^{-n}) - F_{n-1}(x) = -C_{n-1}(y) 2^{-n}\;.
 \end{equation}

Whence, following Lemma \ref{prop:specialslopes} and Case \ref{it:dyadslope} we have for $n\geq J_0$
\BE \label{F-J-Fsoit}
\forall n \geq J_0 \quad
A2^{-\al n} \leq -F(x-2^{-n}) + F(x)\leq B2^{-\al n} \;.
\EE

Let $\rho>0$ and  $J\geq J_0$ such that $2^{-J-1}\leq \rho \leq 2^{-J}$.

Since  $F
 \geq F_J$, then $\Ome_j \subset \Ome$ where $\Ome_J$ is the domain below the graph
of $F_J$. So
\BE  \label{mesureinf}
meas(B(X,\rho) \cap \Ome^c) \leq meas(B(X,\rho) \cap \Ome_J^c)\;.
\EE

But $meas(B(X,\rho) \cap \Ome_J^c)$ is smaller than the area
$h b/2$ of a triangle with altitude $h$ issued from $X$
and a corresponding hypotenuse $b$  (see 
Figures below).


\begin{table}[h] 
\begin{tabular}{cc}  
\includegraphics[width=8cm]{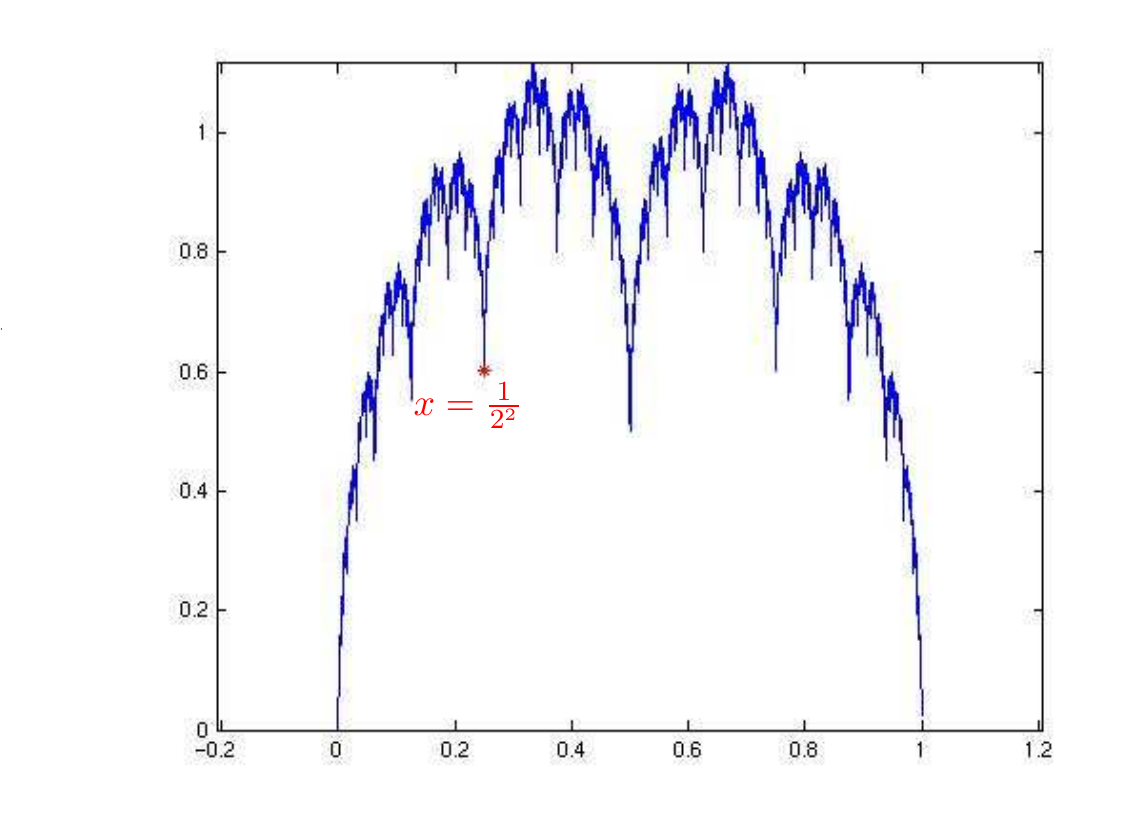}  & \includegraphics[width=8cm]{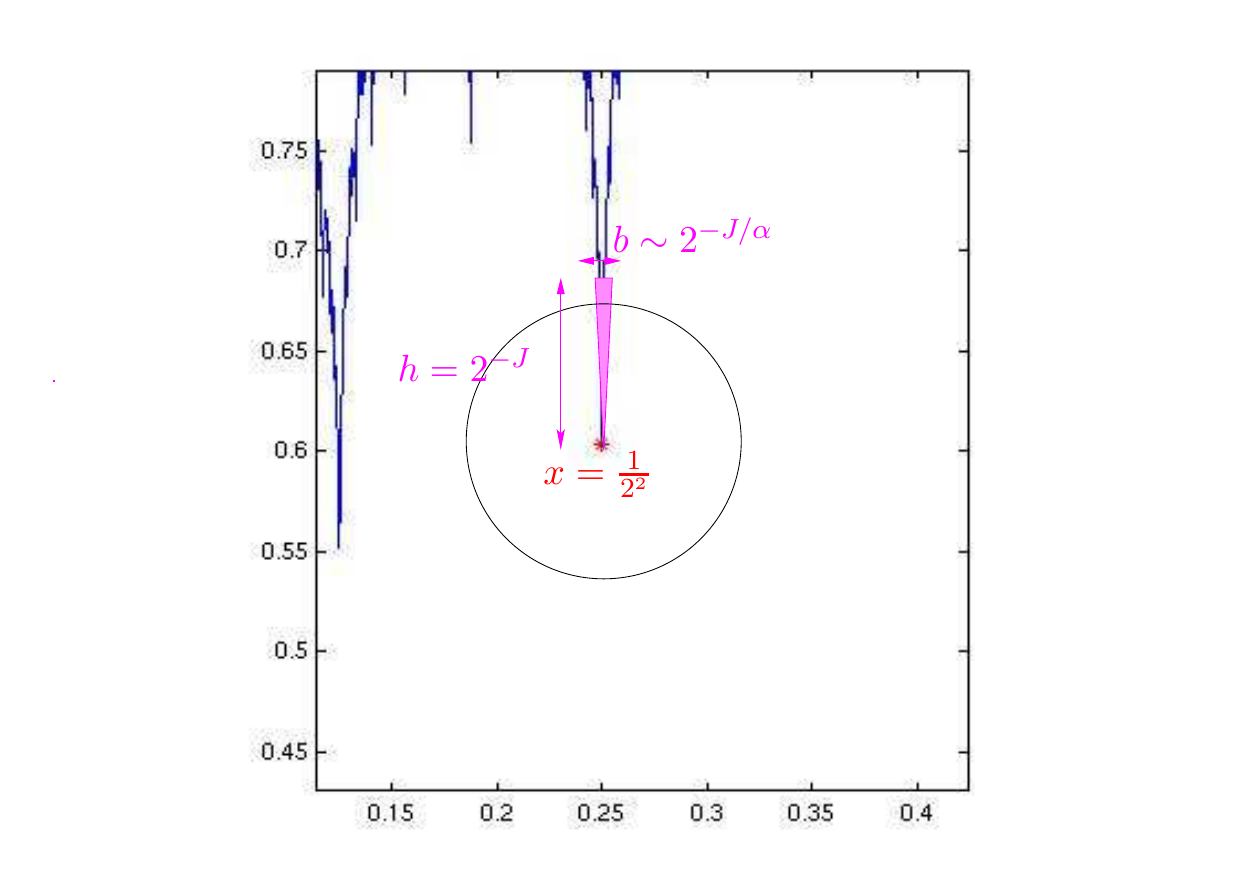}\\
Overview of the function & Zoom at the point $x=\frac{1}{2^2}$\\
\end{tabular} \end{table}

Clearly, we can take $h\sim 2^{-J}$. On the other hand, if we write
$2^{-(j+1)} < b/2 < 2^{-j}$ with $j \geq J$, then
using properties (\ref{F+J-Fsoit}) and (\ref{F-J-Fsoit})
(in which we replace $n\alpha$ by $J$), we get
$b/2 \sim 2^{-J/\al}$. Since $\alpha <1$, Equations (\ref{F+J-Fsoit}) and (\ref{F-J-Fsoit}) are valid with $n=\frac{J}{\alpha}\geq J_0$.

Whence
\BE  \label{mesureinfprecis}
meas(B(X,\rho) \cap \Ome^c) \leq C \rho^{1+\frac{1}{\al}}\;.
\EE
We conclude that
\BE \label{minEw}
E^w_{\Ome^c}(X) \geq \frac{1}{\al} - 1\;.
\EE

Since $E^w_{\Ome^c}(X)\leq E^s_{\Ome^c}(X)\leq \frac{1}{\al} - 1  $ this yields
$$E^s_{\Ome^c}(X)=\frac{1}{\al} - 1  \;.$$

Since $f = 1_\Ome$ then
\BE \label{minup}
u^p_f(X) \geq \ds \frac{1}{p} (\ds \frac{1}{\al} - 1)\;.
\EE

Since $meas(B(X,r))=meas(B(X,r)\bigcap\Omega)+meas(B(X,r)\bigcap\Omega^c)$

we get $E^s_{\Ome}(X)=0$, hence $ E^w_{\Ome}(X)=0$.\\
Whence Proposition \ref{xdyadic}.

\vspace{24pt}
\vspace{24pt}



\subsection{Case of a local maximum of $F$}
We will prove the following proposition.

\begin{prop}\label{xmaxlocal}
Let $X=(x,F(x))$.

If $x$ is a local maximum,
\BE \label{Ewmax}
E^w_{\Ome}(X) = \frac{1}{\al} - 1\mbox{ and }E^w_{\Ome^c}(X) =0
\EE
\BE \label{Esmax}
E^s_{\Ome}(X) = \frac{1}{\al} - 1\mbox{ and }E^s_{\Ome^c}(X) =0
\EE
\BE \label{upmax}
u^p_f(X) = \ds \frac{1}{p} (\ds \frac{1}{\al} - 1)\;.
\EE
\end{prop}

Let $x$ be a local maximum of $F$. There is an interval $I$ containing $x$ such that for all $x'\in I$, $F(x)\geq F(x')$. Let $N$ be such that the dyadic interval $\left[\frac{k_{N}}{2^{N}},\frac{k_N+1}{2^N}\right]$ which contains $x$ is contained in $I$.

Following Lemma \ref{prop:minmax2}, we know that $x$ has the binary expansion (\ref{eq:binexmax}), i.e $x=\frac{k_{N_0}}{2^{N_0}}+\sum\limits_{l=0}^\infty\frac{1}{2^{2l+1+N_0}}$\;.

As a consequence of Lemma \ref{prop:specialslopes}, and following Case (\ref{it:maxslope}), one can find $J_0$ and two constants $A$ and $B$ such that for $n\geq J_0$ Equation (\ref{eq:maxslope}) holds. Remark that it implies clearly that for $n\geq J_0$ $i_{n}(x)=1$ if $n$ is odd, and $i_{n}(x)=0$ if $n$ is even.\\

Our goal now is to evaluate $F(x)-F(x') $ with $x'$ in the interval $\left[\frac{k_n}{2^n},\frac{k_n+1}{2^n}\right] \subset \left[\frac{k_{N}}{2^{N}},\frac{k_{N+1}}{2^{N}}\right]$ and $x' \neq x$. If $x'$ is a dyadic then we take its expansion of type $i_j(x')=0$ for $j$ large enough.\\

Let $m\geq n$ be the smallest integer such that $i_m(x) = i_m(x')$ and $i_{m+1}(x)\neq i_{m+1}(x')$.
To fix the ideas, suppose that $i_{m+1}(x)=1 $ and $i_{m+1}(x')=0$. Thus
\begin{equation}
\begin{split}
\frac{1}{32^{m-1}}\geq x-x'\geq \frac{1}{2^{m+1}}+\frac{1}{2^{m+3}}-\sum\limits_{j=m+2}^\infty\frac{1}{2^j}\geq \frac{1}{2^{m+3}}\;.
\end{split}
\end{equation}

Since $i_{m}(x)=0 $ we have $ C_{m-1}(x)>0$.\\

Thus
\begin{equation}
F(x)-F(x')=\underbrace{C_{m-1}(x)(x-x')}_{(I)}+
\underbrace{\sum\limits_{k=m}^{+\infty}2^{-k\alpha}\left(\Lambda(\tau^kx)-\Lambda(\tau^kx')\right)}_{(II)}\;.
\end{equation}

We have
\begin{equation}
A2^{m(1-\alpha)}2^{-m-3}=C_12^{-m\alpha}\leq (I)\leq B2^{m(1-\alpha)}2^{-m-1}/3=C_22^{-m\alpha}\;.
\end{equation}
Since for $k\geq m-1$ we have $\Lambda(\tau^kx)=1/3 $, this yields
\begin{equation}
0\leq (II)\leq C_32^{-m\alpha}\;.
\end{equation}

Thus we have
\begin{equation}\label{eq:maxexpom}
C_12^{-m\alpha}\leq F(x)-F(x')\leq (C_3+C_2)2^{-m\alpha}\;.
\end{equation}

Let us compute the weak and strong exponents at $x$.\\

Let $\rho$ and $J\geq J_0$ such that $2^{-J-1}\leq \rho\leq 2^{-J}$.
Thus obviously
$$meas(B(X,\rho)\bigcap\Omega)\leq meas(B(X,2^{-j})\bigcap\Omega) \;.$$

Remark first that if $(x',y') \in B(X,2^{-J})\bigcap \Omega $ then $|x-x'| < 2^{-J}$, $|y'-F(x)| < 2^{-J}$
and $y' \leq F(x')$.
Since $x$ is a local maximum on the interval $\left[\frac{k_J}{2^J},\frac{k_J+1}{2^J}\right] $, then $y' \leq F(x') \leq F(x)$
and so $0 \leq F(x)-F(x') < 2^{-J}$. Hence $(x',F(x')) \in B(X,2^{-J})\bigcap \Omega $.

Furthermore since $F(x')$ satisfies $0\leq F(x)-F(x')\leq 2^{-J}$, and following Equation (\ref{eq:maxexpom}) $x'$ belongs to $[x-C2^{-J/\alpha},x+C2^{-J/\alpha} ]$ with $C$ depending only on $C_3+C_2 $.

Thus $B(X,2^{-J})\bigcap \Omega $ is contained in a rectangle of length $2^{-J} $ and width $C2^{-J/\alpha} $.

This yields
\begin{equation}
meas(B(X,\rho)\bigcap\Omega)\leq C2^{-J(1+\frac{1}{\alpha})}\leq C'\rho^{1+\frac{1}{\alpha}}\;.
\end{equation}
We can conclude that

$$E^w_\Omega(X)\geq  \frac{1}{\alpha}-1\;.$$
 Since $E^w_\Omega(X)\leq E^s_\Omega(X)\leq  \frac{1}{\alpha}-1$ this yields

\begin{equation}
E^w_\Omega(X)=E^s_\Omega(X)=\frac{1}{\alpha}-1\;.
\end{equation}

Since $meas(B(X,\rho))= meas(B(X,\rho)\bigcap\Omega)+meas(B(X,\rho)\bigcap\Omega^c)$ we get
$$E^s_{\Omega^c}(X)=E^w_{\Omega^c}(X)=0.$$

And finally
\begin{equation}
u_f^p(X)=\frac{1}{p}\left(\frac{1}{\alpha}-1\right)\;.
\end{equation}




\subsection{Case of $x \notin \mathcal{D} \bigcup \mathcal{S}$}\label{sec:rx}
If $x \notin \mathcal{D} \bigcup \mathcal{S}$ then we will compute separately the weak and strong exponents.
We will first prove that for any point $x$ in $[0,1]$ which is not a maximum or a minimum of $F$ the two weak exponents vanish.\\

\begin{prop}\label{prop:weakexponents}
Let $x \notin \mathcal{D} \bigcup \mathcal{S}$ and $X=(x,f(x))$.\\
Then $E^w_\Omega(X)=E^w_{\Omega^c}(X)=0 $.

\end{prop}

{\bf Proof.}
We will prove first that we have always $E^w_{\Omega^c}(X)=0 $, but will separate the proofs in cases $r(x)>1 $ and $r(x)=1$. Then we will prove that $E^w_\Omega(X)=0 $ and prove it separately for $s(x)>1$, and $s(x)=1$.\\
\begin{itemize}
\item Case $r(x)>1$. We follow the notations of Case \ref{it:r1slope} of Proposition \ref{prop:specialslopes}, i.e one can find two subsequences $J_n$ and $J'_n$ such that $\frac{J'_n}{J_n}>1 $ and $i_{J_n}(x)=i_{J'_n+1}(x) $, $i_j(x)+i_{J_n}(x)=1 $ for $J_n<j<J_n'+1$.
Suppose without loose of generality that $i_{J_n+1}(x)=0 $. Let $\tilde{x}_{n}=\frac{K_{J_n}}{2^{J_n}}=\sum\limits_{j=1}^{J_n}\frac{i_j(x)}{2^j} $. Thus we have
\begin{equation}
2^{-J'_n-1}\leq x-\tilde{x}_n\leq 2^{-J'_n}\;.
\end{equation}

Since Case \ref{it:r1slope} of Proposition \ref{prop:specialslopes} holds, we get

\begin{equation}\label{eq:pentersup1}
\begin{split}
A2^{(1-\alpha)J'_n}2^{-J'_n-1}&\leq F_{J'_n-1}(x)-F_{J'_n-1}(\tilde{x}_n)\leq B2^{(1-\alpha)J'_n}2^{-J'_n}\\
A'2^{-\alpha J'_n}&\leq F_{J'_n-1}(x)-F_{J'_n-1}(\tilde{x}_n)\leq B'2^{-\alpha J'_n}\;.
\end{split}
\end{equation}

We have $F(x)=F_{J'_n-1}(x)+\underbrace{\sum\limits_{k=J'_n}^{+\infty}2^{-k\alpha}\Lambda(\tau^k x)}_{\geq 0}$ and $F_{J'_n-1}(\tilde{x}_n)=F(\tilde{x}_n)$. Thus following (\ref{eq:pentersup1}) we have

\begin{equation}\label{eq:pentersup11}
\begin{split}
A2^{(1-\alpha)J'_n}2^{-J'_n-1}&\leq F(x)-F(\tilde{x}_n)\leq B2^{(1-\alpha)J'_n}2^{-J'_n}+2^{-J'_n\alpha}\sum\limits_{k=0}^{+\infty}2^{-k\alpha}\Lambda(\tau^{k+J'_n} x)\\
A'2^{-\alpha J'_n}&\leq F(x)-F(\tilde{x}_n)\leq B'2^{-J'_n}+2^{-\alpha J'_n}F(\tau^{J'_n}x)\\
A'2^{-\alpha J'_n}&\leq F(x)-F(\tilde{x}_n)\leq B'2^{-J'_n}+2^{-\alpha J'_n}F(\tau^{J'_n}1/3)\\
A'2^{-\alpha J'_n}&\leq F(x)-F(\tilde{x}_n)\leq C2^{-\alpha J'_n}
\end{split}
\end{equation}

indeed the maximum of $F$ is reached at abscissas $1/3$ or $2/3$.\\

We can now apply the mean value theorem and get that for each $n\geq J_0$ we can find $y_n\in ]x-2^{-J'_n},x+2^{-J'_n}[$ such that $F(x)-F(y_n)=A'.2^{-\alpha n}/2$.

Thus using Lemma \ref{lem:Calphaheur1} we can conclude that $E^w_{\Omega^c}(X)=0$.

\item Case $r(x)=1$.\\ Let $J_n$ be defined just as in Lemma \ref{prop:slopes}, i.e that one can find $\delta>0$, $\delta'>0$ and $J_n$ such that Equation (\ref{eq:slopes}) is satisfied.\\

Following the definition of $r(x)$, for all $\gamma>0$ there exists $n_0$ such that for all $j\geq J_{n_0}$ $\left|K_{j}2^{-j}-x\right|>2^{-j(1+\gamma)} $. Thus in particular for all $n\geq n_0$ we have

\begin{equation}
2^{-J_n}>\left|K_{J_n}2^{-J_n}-x\right|>2^{-J_n(1+\gamma)}\;.
\end{equation}

Suppose on one hand $C_{J_n-1}(x)\geq 0 $. Then choose $\tilde{x}_n=K_{J_n}2^{-J_n}$ if $x\in  ]K_{J_n}2^{-J_n},2^{-J_n}+K_{J_n}2^{-J_n}[$ (respectively $\tilde{x}_n=K_{J_n}2^{-J_n}-2^{-J_n}$ if $ x\in  ]K_{J_n}2^{-J_n}-2^{-J_n},K_{J_n}2^{-J_n}[$).

We have obviously
\begin{equation}
F_{J_n-1}(x)-F_{J_n-1}(\tilde{x}_n)=C_{J_n-1}(x)(x-\tilde{x}_n)\geq 0\;.
\end{equation}
If we suppose on the other hand $C_{J_n-1}(x)\leq 0 $, then we can choose in the same way a dyadic number
$\tilde{x}_n=\frac{k}{2^{J_n}}$ such that
\begin{equation}
F_{J_n-1}(x)-F_{J_n-1}(\tilde{x}_n)=C_{J_n-1}(x)(x-\tilde{x}_n)\geq 0\;.
\end{equation}

Together with Equation (\ref{eq:slopes}) this yields in any of these cases that
\begin{equation}
\begin{split}
\delta' 2^{-\alpha J_n}&\geq \left|F_{J_n-1}(x)-F_{J_n-1}(\tilde{x}_n)\right|\geq \delta 2^{-J_n(1+\gamma)}2^{(1-\alpha)J_n}\\
\delta' 2^{-\alpha J_n}&\geq F_{J_n-1}(x)-F_{J_n-1}(\tilde{x}_n)\geq \delta 2^{-J_n(1+\gamma)}2^{(1-\alpha)J_n}\;.
\end{split}
\end{equation}
Since $\sum\limits_{k=J'_n}^{+\infty}2^{-k\alpha}\Lambda(\tau^k x)\geq 0$ and $F_{J_n-1}(\tilde{x}_n)= F(\tilde{x}_n)$ we get

\begin{equation}
\begin{split}
\delta' 2^{-\alpha J_n}+2^{-\alpha J_n}F(1/3)&\geq F(x)-F(\tilde{x}_n)\geq F_{J_n-1}(x)-F_{J_n-1}(\tilde{x}_n)\\
&\geq \delta 2^{-J_n(1+\gamma)}2^{(1-\alpha)J_n}\\
C2^{-\alpha J_n}&\geq  F(x)-F(\tilde{x}_n)\geq\delta 2^{-J_n(\gamma+\alpha)}\;.
\end{split}
\end{equation}

To get  $E^w_{\Omega^c}(X)$ we only have to adapt the proof of Lemma \ref{lem:Calphaheur1} to the case $r_n=2^{-J_n}$. Suppose without lost of generality that $x<\tilde{x}_n$ (the other case can be treated in a similar way) and let $r_n=2^{-J_n}$ for $n\geq n_0$.\\

Indeed, since for $\gamma$ small enough and for $n$ large enough
$2^{-J_n}$ is negligeable in front of $2^{-(\alpha+\gamma)J_n}$ (what we denote
$2^{-(\alpha+\gamma)J_n}>>2^{-J_n}$), following the mean value theorem we can find $b_n\in ]\min(x,\tilde{x}_n),\max(x,\tilde{x}_n)[ $ such that $b_n=\sup\{u_n\in]x,\tilde{x}_n[, f(u_n)=f(x)-r_n\}$. For all $t\in ]b_n,\tilde{x}_n[$, we have $f(t)<f(x)-r_n$. Thus following the same method as in Lemma  \ref{lem:Calphaheur1} we can find $C>0$ such that
\begin{equation}
C'r_n^2\geq meas\left(B(X,r_n)\bigcap\Omega^c\right)\geq Cr_n^{(1+\frac{\gamma}{\alpha})+1}\;.
\end{equation}

This yields
\begin{equation}
2\leq \liminf\limits_{n\rightarrow+\infty}\frac{\log\left( meas\left(B(X,r_n)\bigcap\Omega^c\right)\right)}{\log(r_n)}\leq
\limsup\limits_{n\rightarrow+\infty}
\frac{\log\left( meas\left(B(X,r_n)\bigcap\Omega^c\right)\right)}{\log(r_n)}\leq 2+\frac{\gamma}{\alpha}\;.
\end{equation}

Since $\gamma>0$ is arbitrary and $r_n$ is independent of $\gamma$, we have the result and $E^w_{\Omega^c}(X)=0$.

\item Case $s(x)>1$.\\
Following Case \ref{it:s1slope}, then one can find two subsequences $J_n$ and $J'_n$ with $
\frac{J'_n}{J_n}>1$ for all $n$, such that $i_j(x)+i_{j+1}(x)=1$ for $J_n<j<J'_n$ and $i_{J'_n}(x)=i_{J'_n+1}(x) $.
%
Suppose without loosing generality that $i_{J'_n}(x)=0 $.

Let $\tilde{X}_n $ such that $\tilde{X}_n=\frac{k_{J_n}}{2^{J_n}}+\frac{2}{3(2^{J_n})}=
\sum\limits_{j=1}^{J_n}\frac{i_j(x)}{2^j}+\frac{2}{3(2^{J_n})}$. We have clearly

\begin{equation}
2^{-J'_n-1}\leq -x+\tilde{X}_n\leq 2^{-J'_n}\;.
\end{equation}

Following the same sketch as in the proof with $r(x)>1$ we can say that, using Case 4 of Proposition \ref{prop:slopes}

\begin{equation}
A'2^{-\alpha J'_n}\leq F_{J'_n-1}(\tilde{X}_n)- F_{J'_n-1}(x)\leq B'2^{-\alpha J'_n}
\end{equation}
and since $2^{-\alpha J'_n}F(\tau^{J'_n}1/3)\geq \sum\limits_{k=J'_n}^{+\infty}2^{-k\alpha}\Lambda(\tau^k\tilde{X}_n)- \sum\limits_{k=J'_n}^{+\infty}2^{-k\alpha}\Lambda(\tau^kx)\geq 0 $ we have indeed
\begin{equation}
A'2^{-\alpha J'_n}\leq F(\tilde{X}_n)-F(x)\leq C2^{-\alpha J'_n}\;.
\end{equation}

Thus using the mean value theorem and Lemma \ref{lem:Calphaheur2} as in the previous case we conclude that $E^w_\Omega(X)=0$.

\item Case $s(x)=1$.

Let $J_n$ be defined just as in Lemma \ref{prop:slopes}, i.e that one can find $\delta>0$, $\delta'>0$ and $J_n$ such that Equation (\ref{eq:slopes}) is satisfied as well as Point \ref{it:slope2} of Lemma \ref{prop:slopes}.\\

Following the definition of $s(x)$, for all $\gamma>0$ there exists $n_0$ such that for all $j\geq J_{n_0}$ $\left|m_j-x\right|>2^{-j(1+\gamma)} $. Thus in particular for all $n\geq n_0$ we have

\begin{equation}
2^{-J_n}>\left|m_{J_n}-x\right|>2^{-J_n(1+\gamma)}\;.
\end{equation}

Suppose on one hand that $C_{J_n-1}(x)> 0 $. Then take $\tilde{X}_n=\frac{K_{J_n}}{2^{J_n}}+\frac{2}{3(2^{J_n})}$. Since $i_{J_n+1}(x)=0 $ we have $F_{J_n-1}(\tilde{X_n})-F_{J_n-1}(x)>0$.


Together with Equation (\ref{eq:slopes}) this yields that
\begin{equation}
\begin{split}
\delta' 2^{-\alpha J_n}&\geq \left|F_{J_n-1}(x)-F_{J_n-1}(\tilde{X}_n)\right|\geq \delta 2^{-J_n(1+\gamma)}2^{(1-\alpha)J_n}\\
\delta' 2^{-\alpha J_n}&\geq -F_{J_n-1}(x)+F_{J_n-1}(\tilde{X}_n)\geq \delta 2^{-J_n(1+\gamma)}2^{(1-\alpha)J_n}\;.
\end{split}
\end{equation}
The same computation as previously yields
\begin{equation}
\begin{split}
\delta' 2^{-\alpha J_n}+2^{-\alpha J_n}F(1/3)&\geq -F(x)+F(\tilde{X}_n)\geq -F_{J_n-1}(x)+F_{J_n-1}(\tilde{X}_n)\\
&\geq \delta 2^{-J_n(1+\gamma)}2^{(1-\alpha)J_n}\\
C2^{-\alpha J_n}&\geq  -F(x)+F(\tilde{X}_n)\geq\delta 2^{-J_n(\gamma+\alpha)}\;.
\end{split}
\end{equation}

To get  $E^w_{\Omega}(X)$ we only have to adapt the proof of Lemma \ref{lem:Calphaheur2}  in the same way we adapt the one of Lemma \ref{lem:Calphaheur1} in the case $r(x)=1$.\\


Thus taking $r_n=2^{-J_n}$ and following the same method as previously we can find $C>0$ such that
\begin{equation}
C'r_n^2\geq meas\left(B(X,r_n)\bigcap\Omega\right)\geq C2^{-J_n(1+\frac{\gamma}{\alpha})}2^{-J_n}\;.
\end{equation}
This yields
\begin{equation}
2\leq \liminf\limits_{n\rightarrow+\infty}
\frac{\log\left( meas\left(B(X,r_n)\bigcap\Omega\right)\right)}{\log(r_n)}\leq
\limsup\limits_{n\rightarrow+\infty}\frac{\log\left( meas\left(B(X,r_n)\bigcap\Omega\right)\right)}{\log(r_n)}
\leq 2+\frac{\gamma}{\alpha}\;.
\end{equation}

Since $\gamma>0$ is arbitrary, we have the result and $E^w_{\Omega}(X)=0$.

Hence the proof of Proposition \ref{prop:weakexponents}.

\end{itemize}

For what concerns the strong accessibility exponent we have the following result.

\begin{prop}\label{prop:strongexponents}
Suppose $x \notin \mathcal{D} \bigcup \mathcal{S}$  and let $X_0=(x,F(x))$. Then
\begin{enumerate}
\item\label{it:stromc} If $r(x)>\frac{1}{\alpha}$ then $E^s_{\Omega^c}(X_0)=\frac{1}{\alpha}-1 $.
\item If $s(x)>\frac{1}{\alpha}$ then $E^s_{\Omega}(X_0)=\frac{1}{\alpha}-1 $.
\item Let $D_\alpha$ the set of $x \notin \mathcal{D} \bigcup \mathcal{S}$ such that $r(x)>\frac{1}{\alpha} $ and $s(x)>\frac{1}{\alpha} $. Then the Hausdorff dimension of $D_\alpha$ is $\alpha$.
\end{enumerate}

\end{prop}

{\bf Proof.}
\begin{enumerate}
\item Let us prove Point \ref{it:stromc}. Since $r(x)>\frac{1}{\alpha}$, and following Point \ref{it:r1slope} of Lemma \ref{prop:specialslopes}, for $\delta>0$ such that $r(x)-\delta>\frac{1}{\alpha} $ we can find $J'_n$ and $J_n$ such that
\begin{itemize}
\item $x_n=K_{J_n}2^{-J_n}$ and $|x-x_n|\leq 2^{-J'_n} $.
\item $|F(x)-F(x_n)|\leq 2^{-\alpha J'_n} $.
\end{itemize}
Since we can choose $J'_n$ such that  $2^{-J'_n}\leq C2^{(r(x)-\delta)J_n}$
(see the proof of Point \ref{it:r1slope} of Lemma \ref{prop:specialslopes}), then
$|x-x_n|$ is negligeable in front of $2^{-J_n}$ (what we denote $|x-x_n|<<2^{-J_n} $)
and $|F(x)-F(x_n)|\leq 2^{-\alpha J'_n}<<2^{-J_n} $.\\

Thus we can choose a constant $C$ such that with $\rho_n=C2^{-J_n}$ and $B(X_0,\rho_n/2)\subset B((x_n,F(x_n)),\rho_n) $.
Following the proof of Proposition \ref{xdyadic} and more precisely Equation (\ref{mesureinfprecis}) we have
\begin{eqnarray*}
meas((B(X,\rho_n/2)\bigcap \Omega^c)&\leq meas (B((K_{J_n}2^{-J_n},F(K_{J_n}2^{-J_n}),\rho_n)\bigcap \Omega^c))\\
&\leq C\rho_n^{1+1/\alpha}\;.
\end{eqnarray*}

This yields $E^s_{\Omega^c}(X)\geq \frac{1}{\alpha}-1$.

Since $E^s_{\Omega^c}(X)\leq \frac{1}{\alpha}-1$, we get $E^s_{\Omega^c}(X)=\frac{1}{\alpha}-1$.\\
\item We follow exactly the same proof as previously replacing $x_n$ by $\tilde{X_n}$ the sequence of local maxima defined in the proof of Point \ref{it:s1slope} of  Lemma \ref{prop:specialslopes}.

\item We follow here the results proved by \cite{Dur} and summarized in \cite{Bug} for our special case.
Indeed recall the definition given in \cite{Dur} of an ubiquitous system in  a real interval of $\R$.\\

\begin{defi}
Let $U$ be a real open interval. Let $(x_i)_{i\geq 1}$ be points in $U$ and let $(r_i)_{i\geq 1}$ be a sequence of positive real numbers such that $\lim\limits_{i\rightarrow\infty}r_i=0$. The family $(x_i,r_i)_{i\geq 1} $ is a homogeneous ubiquitous system in $U$ if the set $\limsup\limits_i B(x_i,r_i) $ is of full Lebesgue measure in $U$.
\end{defi}

Theorem D of \cite{Bug} proved in \cite{Dur} yields the following result.
\begin{thm}\label{th:bug}
Let $\tau$ be a real number with $\tau \geq 1$. With the above notations if the families $(x_i,r_i)_{i\geq 1} $ and $(x'_i,r'_i)_{i\geq 1} $ are two homogeneous ubiquitous systems in U, then the Hausdorff dimension of the set $\limsup B(x_i,r_i^\tau)\bigcap \limsup B(\tilde{x}_j,\tilde{r}_j^\tau)$ is at least equal to $\frac{1}{\tau}$.
\end{thm}

Let $U=]0,1[$ and consider $\mathcal{K}_1=\{(\frac{k}{2^j},2^{-j}),k\in\N,0< k<2^j,j\geq 1\}$. It is a countable set and can be written as $\mathcal{K}_1=\{(x_i,r_i),i\geq 1\} $ with $x_i$ a dyadic number for all $i\geq 1$. Let $\mathcal{K}_2=\{(x,r),x\in\mathcal{S},r=\frac{2^{-j}}{3}\mbox{ for }j\geq 1\}$. It is again a countable set and we can rewrite it as $\mathcal{K}_2=\{(\tilde{x}_i,\tilde{r}_i),i\geq 1\} $ with $\tilde{x}_i\in\mathcal{S}$ for all $i\geq 1$.\\
It is clear that $\limsup B(x_i,r_i)$ and $\limsup B(\tilde{x}_i,\tilde{r}_i)$ are of full Lebesgue measure.\\

Remark then that $D_\alpha=\limsup B(x_i,r_i^\tau)\bigcap \limsup B(\tilde{x}_i,\tilde{r}_i^\tau) $
with $\tau=\frac{1}{\alpha}$. Since $D_\alpha\subset \bigcup\limits_{j\geq J,0\leq k\leq 2^j}
B(x_i,r_i^{\tau})$ the Hausdorff dimension of $D_\alpha$ is less or equal than $\alpha $.
We apply Theorem \ref{th:bug} and we find it is exactly $\alpha$.\\

\end{enumerate}

\vspace{24pt}
\subsection{Proof of Theorem \ref{th:results}}
Propositions \ref{xdyadic}, \ref{xmaxlocal},
\ref{prop:weakexponents} and \ref{prop:strongexponents} achieve  the proof of Theorem \ref{th:results}.

\vspace{24pt}

\subsection*{Acknowledgements}
M. Ben Slimane was supported by the Research Center, College of Science, King Saud University.

\end{document}